\documentclass[reqno, 11pt]{article}
\usepackage{amssymb,dsfont,graphicx,bbm}
\usepackage{amsthm,amsmath}
\usepackage[usenames,dvipsnames]{color}
\usepackage{fullpage}
\usepackage{ifpdf}
\usepackage[pdftex]{thumbpdf}       
\graphicspath{{graphs/}}

\ifpdf
\usepackage[pdftex]{hyperref}
\hypersetup{%
pdftitle={A two station queue with dependent preparation and service times},
pdfauthor={I.J.B.F. Adan, O.J. Boxma, M. Vlasiou},
pdfkeywords={},
bookmarks=true,%
bookmarksnumbered=true,
pdfstartview={FitV},
pdfpagelayout={TwoColumnLeft},
colorlinks=true,
citecolor=Plum,
linkcolor=RedOrange,
urlcolor=RawSienna,
bookmarksopen=true,%
bookmarksopenlevel=1,     
unicode=true,%
breaklinks=true,%
}%
\else
\newcommand\phantomsection\relax
\newcommand{\url}[1]{#1}
\newcommand{\href}[2]{#2}
\usepackage{nohyperref}
\fi

\theoremstyle{plain}              
\newtheorem{theorem}{Theorem}

\theoremstyle{remark}
\newtheorem{remark}{Remark}

\numberwithin{equation}{section}    

\newcommand{\m}[1]{\mathcal{#1}}
\newcommand{\e}{\mathbb{E}}
\newcommand{\p}{\mathbb{P}}

\newcommand{\cov}{\ensuremath{\mathrm{cov}}}

\newcommand{\Dfb}[1][]{\ensuremath{F_B^{#1}}} 
\newcommand{\Dfa}[1][]{\ensuremath{F_A^{#1}}}
\newcommand{\Dfw}[1][]{\ensuremath{F_W^{#1}}} 

\newcommand{\dfw}[1][]{\ensuremath{f_W^{#1}}}

\newcommand{\lta}[1][]{\ensuremath{\alpha^{#1}}}
\newcommand{\ltw}[1][]{\ensuremath{\omega^{#1}}}

\begin{document}
\title{A two-station queue with dependent preparation and service times}
\author{M.\ Vlasiou$^{*}$, I.J.B.F.\ Adan$^{**,***}$, and O.J.\ Boxma$^{**,***}$}
\date{September 11, 2007}
\maketitle

\begin{center}
$^{*}$ Georgia Institute of Technology,
\\H.\ Milton Stewart School of Industrial \& Systems Engineering,
\\765 Ferst Drive, Atlanta GA 30332-0205, USA.
\end{center}

\begin{center}
$^{**}$ Eindhoven University of Technology,
\\Department of Mathematics \& Computer Science,
\\P.O.\ Box 513, 5600 MB Eindhoven, The Netherlands.
\end{center}

\begin{center}
$^{***}$ EURANDOM,\\
P.O.\ Box 513, 5600 MB Eindhoven, The Netherlands.
\\\vspace{0.3cm}\href{mailto:vlasiou@gatech.edu}{vlasiou@gatech.edu}, \href{mailto:iadan@win.tue.nl}{iadan@win.tue.nl}, \href{mailto:boxma@win.tue.nl}{boxma@win.tue.nl}
\end{center}

\begin{abstract}
We discuss a single-server multi-station alternating queue where the preparation times and the service
times are auto- and cross-correlated. We examine two cases. In the first case, preparation and service times depend
on a common discrete time Markov chain. In the second case, we assume that the service times depend on the
previous preparation time through their joint Laplace transform. The waiting time process is directly analysed by
solving a Lindley-type equation via transform methods. Numerical examples are included to demonstrate the effect
of the autocorrelation of and the cross-correlation between the preparation and service times.
\end{abstract}

\section{Introduction}\label{s:intro}
Consider a server alternating between two service points, serving one customer at a time. At each service point there is an infinite queue of customers waiting for service. Before being served by the server, a customer first must undergo a preparation phase, which starts immediately after the server has completed service at that particular service point and has moved to the next one. The server is not involved at all with the preparation phase. The server is obliged to alternate; therefore he serves all odd-numbered customers at one service point and all even-numbered customers at the other. Thus the server, after having finished serving a customer at one service point, may have to wait for the preparation phase of the customer at the other service point to be completed. Let $W_n$ be the time the server has to wait before he can start serving the $n$-th customer. If $B_n$  is the preparation time of the $n$-th customer and $A_n$ is the service time of the $n$-th customer, then  $W_n$ can be defined recursively by
\begin{equation}\label{recursion}
W_{n+1}=\max\{0, B_{n+1}-A_n-W_n\}, \quad n\geqslant 1.
\end{equation}

This alternating service model occurs in many applications. For example, this strategy is followed by surgeons performing eye surgeries. Another example where \eqref{recursion} occurs comes from inventory theory, in particular from the analysis of two-carousel systems. This application is considered in \cite{hassini03,koenigsberg86,park03}. Equation~\eqref{recursion} is introduced in \cite{park03} and has been studied further in \cite{vlasiou07a,vlasiou05,vlasiou05b,vlasiou04,vlasiou07}; see these references for a further discussion on applications. A restriction of the previous work is that $\{A_n\}$ and $\{B_n\}$ are two mutually \emph{independent} sequences of i.i.d.\ random variables. In the present study, we analyse the above queuing model under two different dependence structures.

In the first dependence structure we study, the distributions of the preparation and service times are regulated by an irreducible discrete-time Markov chain. Specifically, we assume that each transition of the Markov chain generates a new preparation time and its corresponding service time. Given the state of the Markov chain at times $n$ and $n+1$, the distributions of $A_n$ and $B_{n+1}$ are independent of one another for all $n$. However, the distributions of $A_n$ and $B_n$ depend on the state of the Markov chain. Such a dependence structure occurs naturally in many applications. For example, in the application involving two carousels that is described in \cite{park03,vlasiou05,vlasiou04}, one can intuitively see that if an order consists of multiple items stored on one carousel, then there are strategies for the preparation of the carousel, where a long preparation time $B_n$ implies that the service time $A_n$ will be relatively short, while being independent of all other past or future preparation and service times. Note that the service time $A_n$ is defined as the time necessary to pick an order on that carousel, which is comprised of the time needed to pick each individual item in one order plus the time needed to rotate between the items of one order. Also in this example, the service and preparation times may depend on the carousel. That is, if one carousel is slower than the other, then the service and preparation times on the slower carousel will be longer than on the faster one. Hence, the service and preparation times are strictly periodic, which is a special case of Markov-dependence. Also queueing models with strictly periodic arrivals have been frequently studied; see, \cite{cohen95,lemoine89,rolski87,rolski89}. These models arise, for example, in the modelling of inventory systems using periodic ordering policies; see \cite{smits04} and \cite{wagner04}.

The second dependence structure we study assumes that the random variables $A_n$ and $B_{n+1}$ have a joint distribution. In particular, given the length of the service time $A_n$, the following preparation time has a Laplace-Stieltjes transform of a specific form. The form we choose is rather general and allows for various specific dependence structures and preparation time distributions. Later on, we give a few specific examples. Dependencies between a service time and the {\em following} preparation time are also possible in applications. Again for the carousel model described in \cite{vlasiou05,vlasiou04} with orders consisting of multiple items, one can have the situation where a ``smart'' preparation strategy is followed, which anticipates the expected delay of the server for the previous order. Thus, knowing that the previous service time is relatively long, the other carousel rotates at a starting point that may be further away, but reduces the service time of the following order.

One should note here that \eqref{recursion} is, up to the minus sign in front of $W_n$, equal to Lindley's recursion, which is one of the most important and well-studied recursions in applied probability. Both dependence structures studied are also motivated by analogous cases studied for Lindley's recursion. In Lindley's recursion, $A_n$ represents the interarrival time between customers $n$ and $n+1$, $B_n$ is the service time of the $n$-th customer, and $W_n$ represents the time the $n$-th customer has to wait before starting his service; see, e.g., \cite{asmussen-APQ} and \cite{cohen-SSQ} for a comprehensive description. We shall compare Lindley's recursion to \eqref{recursion} for the two dependence structures mentioned above.

Queuing models with dependencies between interarrival and service time have been studied by several authors. A review of the early literature can be found in Bhat~\cite{bhat69}. Such dependencies arise naturally in various applications. For example, the phenomenon of dependence among the interarrival times in the packet streams of voice and data traffic is well known; see, e.g., \cite{heffes80,heffes86,sriram86}. However, in \cite{fendick89} the authors argue that in packet communication networks one should also expect two additional forms of dependence: between successive service times and among interarrival times and service times. These forms of dependence occur because of the presence of bursty arrivals and multiple sources with different mean service times (due to different packet lengths), and they may have a dominant effect on waiting times and queue lengths. In this paper, we study such forms of dependence for the model described by \eqref{recursion}.

The presentation is organised as follows. In Section~\ref{7s:Markov-modulated dependencies} we derive the steady-state waiting-time distribution $\Dfw$ for the first case studied. We assume that the service times follow some general distribution that depends on the state of the Markov chain. For the preparation times, in Section~\ref{7ss:Exponential preparation times} we assume that they are exponentially distributed (with a rate depending on the state of the Markov chain), while in Section~\ref{7ss:Phase-type preparation times} we extend the analysis to mixed-Erlang distributions. We complement the results with various numerical examples in Section~\ref{ss:egs1}. In Section~\ref{7s:B_n+1 depends on A_n} we derive $\Dfw$ for the second case studied. We assume that the  service times are exponentially distributed, although, as we remark later on, the analysis can be extended to mixed-Erlang distributions.

\section{Markov-modulated dependencies}\label{7s:Markov-modulated dependencies}
In this section we define the setting and derive the autocorrelation function of preparation or service times, as well as the crosscorrelation of these. We study the case where the preparation times and the service times depend on a common discrete-time Markov chain. This model allows dependencies between preparation and service times. The waiting time in this case is directly derived by using Laplace transforms. For Lindley's recursion, the analogous model has been analysed in Adan and Kulkarni~\cite{adan03}.

We first introduce some notation. For a random variable $Y$ we denote its distribution by $F_Y$ and its density by $f_Y$. A second index is used when we need to distinguish between the different states of the underlying Markov chain. For example, we denote by $F_{W,j} (x)$ the steady-state probability that the waiting time $W \le x$ and the Markov chain is in state $j$. Also, given an event $E$ we use the convention that $\p[Y\leqslant x\, ; E]=\e[\mathbbm{1}_{[Y\leqslant x]} \cdot \mathbbm{1}_{[E]}]$, and likewise for expectations.

As mentioned before, we assume that the sequences $\{A_n\}$ and $\{B_n\}$ are both autocorrelated and cross-correlated. The nature of the dependence we study in this first part is described below.

The distributions of the preparation and service times are regulated by an irreducible, aperiodic, discrete-time Markov chain $\{Z_n\}$, $n \geqslant 1$, with state space $\{1,2,\ldots, M\}$ and transition probability matrix $\mathbf{P}=(p_{i,j})$. To be exact, we have that
\begin{align}
\nonumber \p[A_n&\leqslant x\,; B_{n+1}\leqslant y\,; Z_{n+1}=j \mid Z_n=i\,; B_n\,; (A_\ell, B_\ell, Z_\ell), 1\leqslant\ell\leqslant n-1]\\
\nonumber &\quad =\p[A_1\leqslant x\,; B_{2}\leqslant y\,; Z_{2}=j \mid Z_1=i]\\
\nonumber &\quad =p_{i,j} \p[A_1\leqslant x\,; B_{2}\leqslant y \mid Z_1=i\,; Z_{2}=j]\\
\label{7eq:2} &\quad =p_{i,j} F_{A\mid i}^{}(x) F^{}_{B\mid j}(y),
\end{align}
where $x,y \geqslant 0$ and where $i,j=1,2,\ldots,M$. Thus, given $Z_n$ and $Z_{n+1}$, the distributions of $A_n$ and $B_{n+1}$ are independent of one another for all $n$.

Before introducing the specific assumptions on the distributions of $A_n$ and $B_n$, we first discuss the stability properties of the system. We assume that there are an $i$ and $j$ such that $\p[X_n < 0, Z_{n+1} = j | Z_n = i] > 0$, where $X_n = B_{n+1} - A_n$. This implies that the Markov chain $(W_n, Z_n)$ is stable. To see that, first notice that a regeneration point occurs \index{regenerative process} if $W_{n+1}=0$ and $Z_{n+1}=j$.
Since the Markov chain $\{Z_n\}$ is finite and irreducible, the time between two occurrences of state $j$ is finite in expectation. For each time state $j$ is reached there is a positive probability that $X_n<0$, which implies that $W_{n+1}=0$. Thus, we will have that in a geometric number of steps both events $W_{n+1}=0$ and $Z_{n+1}=j$ will happen at the same step. In particular, the number of steps is finite in expectation. In other words, we will have a regeneration point within a finite expected time.
Hence, since $\{Z_n\}$ is aperiodic,
the limiting distribution of the Markov chain $(W_n, Z_n)$ exists, and thus we can define for $\mathrm{Re}(s)\geqslant 0$, $n\geqslant 1$, and $j=1,2,\ldots,M$ the transforms
$$
\ltw_{j}(s)=\lim_{n\to\infty}\ltw[n]_{j}(s),
$$
where
$$
\ltw[n]_{j}(s)=\e[\mathrm{e}^{-s W_{n}}\,;Z_{n}=j].
$$

Denote by $\lambda_i^{-1}$ the mean and by $s_i$ the second moment of the service time distribution $F_{A\mid i}^{}$. Analogously, define $\mu_i^{-1}$ as the mean of $F_{B\mid i}^{}$ and $\sigma_i$ as its second moment. Moreover, denote by $\mbox{{\boldmath $\varpi$}}=( \varpi_1, \varpi_2,\ldots, \varpi_M)$ the stationary distribution of the Markov chain $\{Z_n\}$. Then, in steady state, the autocorrelation between $A_m$ and $A_{m+n}$ is given by
\begin{equation}\label{eq:autoA}
\rho[A_m,A_{m+n}]=\rho[A_1,A_{n+1}]=\frac{\sum_{i=1}^M\sum_{j=1}^M \varpi_i \bigl(p_{i,j}^{(n)}-\varpi_j\bigr)\lambda_i^{-1}\lambda_j^{-1}}{\sum_{i=1}^M \varpi_i s_i-\Bigl(\sum_{i=1}^M \varpi_i \lambda_i^{-1}\Bigr)^2},
\end{equation}
where
$$
p_{i,j}^{(n)}=\p[Z_{n+1}=j \mid Z_1=i],\qquad n\geqslant0,\quad 1\leqslant i,j\leqslant M.
$$
A similar expression can be derived for the autocorrelation between preparation times. Since $\mathbf{P}$ is aperiodic, we have that $p_{i,j}^{(n)}$ converges to $\varpi_j$ geometrically as $n$ tends to infinity. In other words, the autocorrelation function approaches zero geometrically fast as the lag goes to infinity. For the cross-correlation between $A_n$ and $B_n$ we have that
$$
\rho[A_n,B_{n}]=\rho[A_1,B_{1}]=\frac{\sum_{i=1}^M\varpi_i \lambda_i^{-1}\mu_i^{-1}-\widehat{\mu}\,\widehat{\lambda}}{\Bigl(\sum_{i=1}^M \varpi_i s_i- \widehat{\lambda}^2 \Bigr)^{1/2}\Bigl(\sum_{j=1}^M\varpi_j \sigma_i- \widehat{\mu}^2\Bigr)^{1/2}},
$$
where $\widehat{\lambda}=\sum_{i=1}^M\varpi_i \lambda_i^{-1}$ and $\widehat{\mu}=\sum_{i=1}^M\varpi_i \mu_i^{-1}$.

In the remainder of this section, we assume that the random variables $A_n$ follow an arbitrary distribution that is independent of the past, given $Z_n$, while $B_n$ follows in general a phase-type distribution that is depending on the state of $Z_n$.

When $\{Z_n\}$ is in state $j$, we denote by $\lta_j$ the Laplace-Stieltjes transform of the service distribution $F_{A\mid j}$. We denote the derivative of order $i$ of a function $f$ by $f^{(i)}$ and we have by definition that $f^{(0)}=f$. For simplicity, in the following section we first derive the steady-state waiting-time distribution for this model in case the preparation times are exponentially distributed. Later on, in Section~\ref{7ss:Phase-type preparation times}, we generalise this result to phase-type preparation times.

\subsection{Exponential preparation times}\label{7ss:Exponential preparation times}
In this section we assume that $F^{}_{B\mid j}(x)=1-\mathrm{e}^{-\mu_j x}$. Note that the rates $\mu_j$ need not be distinct; see also Section~\ref{ss:egs1} for some examples. We are interested in the steady-state waiting-time distribution. The next theorem gives the equations satisfied by the waiting-time densities for every state of $\{Z_n\}$.

\begin{theorem}[Exponential preparation times]\label{7th:MC-exp B}
Let $F^{}_{B\mid j}(x)=1-\mathrm{e}^{-\mu_j x}$. Then, for all $j = 1,\ldots, M$,
\[
F_{W,j} (0) = \varpi_j- c_j,
\]
where $c_j=\sum_{i=1}^M p_{i,j}\,\ltw_{i}(\mu_j)\,\lta_i(\mu_j)$, and the waiting-time density is given by
$$
f_{W,j}^{}(x)=\mu_j c_j\,\mathrm{e}^{-\mu_j x}.
$$
The $M^2$ unknown constants $\ltw_{i}(\mu_j)$ are the unique solution to the following system of linear equations
$$
\ltw_{j}(\mu_\ell)=\varpi_j-\frac{\mu_\ell}{\mu_j+\mu_\ell} \sum_{i=1}^M p_{i,j}\,\ltw_{i}(\mu_j)\,\lta_i(\mu_j),\qquad j,\ell=1,\ldots,M.
$$
\end{theorem}
\begin{proof}
From \eqref{recursion} we obtain the following equation for the transforms $\ltw[n+1]_{j}$, $j=1,\ldots,M$.
\begin{align}\label{7eq:ltw e3is}
\nonumber \ltw[n+1]_{j}(s)&=\e[\mathrm{e}^{-s W_{n+1}}\,;Z_{n+1}=j]=\sum_{i=1}^M\p[Z_n=i]\e[\mathrm{e}^{-s\max\{0,B_{n+1}-A_n-W_n\}}\,;Z_{n+1}=j\mid Z_n=i]\\
\nonumber &=\sum_{i=1}^M\p[Z_n=i] p_{i,j} \biggl(\e[\int_0^{A_n+W_n} f_{B_{n+1}}^{}(x)\,\mathrm{d}x \mid Z_n=i\,;Z_{n+1}=j]+\\
&\quad+\e[\int_{A_n+W_n}^\infty \mathrm{e}^{-s(x-A_n-W_n)} f_{B_{n+1}}^{}(x)\,\mathrm{d}x \mid Z_n=i\,;Z_{n+1}=j]\biggr).
\end{align}
Since $Z_{n+1}=j$, we have that $B_{n+1}$ is now exponentially distributed with rate $\mu_j$. Thus, the above equation becomes
\begin{align*}
\nonumber \ltw[n+1]_{j}(s)&=\sum_{i=1}^M\p[Z_n=i] p_{i,j} \biggl(\e[\int_0^{A_n+W_n} \mu_j \mathrm{e}^{-\mu_j x}\,\mathrm{d}x \mid Z_n=i]+\\
\nonumber &\quad+\e[\int_{A_n+W_n}^\infty \mathrm{e}^{-s(x-A_n-W_n)} \mu_j \mathrm{e}^{-\mu_j x}\,\mathrm{d}x \mid Z_n=i]\biggr)\\
\nonumber &=\sum_{i=1}^M \p[Z_n=i] p_{i,j}\e[1-\mathrm{e}^{-\mu_j(A_n+W_n)}+\frac{\mu_j}{\mu_j+s}\,\mathrm{e}^{-\mu_j(A_n+W_n)} \mid Z_n=i]\\
\nonumber &=\sum_{i=1}^M \p[Z_n=i] p_{i,j} \Bigl(1-\frac{s}{\mu_j+s}\,\e[\mathrm{e}^{-\mu_j(A_n+W_n)} \mid Z_n=i]\Bigr)\\
 &=\sum_{i=1}^M p_{i,j}\Bigl(\p[Z_n=i]-\frac{s}{\mu_j+s}\,\ltw[n]_{i}(\mu_j)\,\lta_i(\mu_j)\Bigr).
\end{align*}
So for $n\to\infty$ we have that $\ltw_j(s)$ is given by
\begin{equation}\label{7eq:ltw B exp}
\ltw_{j}(s)= \varpi_j -\sum_{i=1}^M p_{i,j}\,\ltw_{i}(\mu_j)\,\lta_i(\mu_j)+ \frac{\mu_j}{\mu_j+s}\sum_{i=1}^M p_{i,j}\,\ltw_{i}(\mu_j)\,\lta_i(\mu_j).
\end{equation}
Define the constants
$$
c_j=\sum_{i=1}^M p_{i,j}\,\ltw_{i}(\mu_j)\,\lta_i(\mu_j).
$$
Inverting the Laplace transform $\ltw_j$ yields that the density of the waiting time is given by
$$
f_{W,j}^{}(x)=\mu_j c_j\,\mathrm{e}^{-\mu_j x},
$$
and the corresponding distribution has mass $F_{W,j} (0) =\varpi_j- c_j$ at the origin. For $i,j=1,\ldots,M$, the $M^2$ unknown constants $\ltw_{i}(\mu_j)$ that are needed in order to determine the unknown constants $c_j$ are the unique solution to the system of linear equations given by the expression
\begin{equation}\label{7eq:system}\index{balance equations}
\ltw_{j}(\mu_\ell)=\varpi_j-\frac{\mu_\ell}{\mu_j+\mu_\ell} \sum_{i=1}^M p_{i,j}\,\ltw_{i}(\mu_j)\,\lta_i(\mu_j),\qquad j,\ell=1,\ldots,M;
\end{equation}
see Equation~\eqref{7eq:ltw B exp}. The uniqueness of the solution follows from the general theory of Markov chains, which states that there is a unique stationary distribution and thus also a unique solution to the system of equations formed by \eqref{7eq:system} for all $j,\ell=1,\ldots,M$.
\end{proof}

The result given in Theorem~\ref{7th:MC-exp B} is expected. Evidently, since the preparation time $B_n$ is exponentially distributed (with a rate depending on the state of the Markov chain) and $W_n$ is the residual preparation time, we have that for every state $j$ of the Markov chain, the waiting-time distribution has mass at zero and the conditional waiting time is exponentially distributed with rate $\mu_j$.

Observe that Theorem~\ref{7th:MC-exp B} generalises the statement of Theorem~1 in \cite{vlasiou05}, which gives the steady-state waiting-time density in case $\{A_n\}$ and $\{B_n\}$ are mutually independent sequences of i.i.d.\ random variables and $B$ follows an Erlang distribution. Specifically, if the Markov chain in Theorem~\ref{7th:MC-exp B} has only one state (and thus there is a unique service-time distribution and a unique rate $\mu$ for the exponentially distributed preparation times) and the Erlang distribution $\Dfb$ in \cite[Theorem~1]{vlasiou05} has only one phase, then the statements of these two theorems are identical. Observe, for example, that \eqref{7eq:system} reduces to \cite[Eq.\ (3.2)]{vlasiou05} as now $\varpi_j=1$, and $p_{i,j}=1$.

\subsection{Phase-type preparation times}\label{7ss:Phase-type preparation times}

Assume now that if the Markov chain is in state $j$, the preparation time is with probability $\kappa_n$ equal to a random variable $Y_n$, $n=1,\ldots,N$, that follows an Erlang distribution with parameter $\mu_j$ and $n$ phases. In other words the distribution function of $B$ is given by
\begin{equation}\label{7eq:distr}
F^{}_{B\mid j}(x)=\sum_{n=1}^N \kappa_n\biggl(1-\mathrm{e}^{-\mu_j x}\sum_{\ell=0}^{n-1}\frac{(\mu_j x)^\ell}{\ell!}\biggr) ,\qquad x \geqslant 0.
\end{equation}
This class of phase-type distributions may be used to approximate any given distribution on $[0,\infty)$ for the preparation times arbitrarily close; see Schassberger~\cite{schassberger-W}. The waiting-time density for this case is given by the following theorem.
\begin{theorem}[Mixed-Erlang preparation times]\label{7th:MC-PH B}
Let $F^{}_{B\mid j}$
be given by \eqref{7eq:distr}. Then, for $j = 1, \ldots, M$,
\begin{equation}\label{mass}
F_{W,j} (0) =
\varpi_j-\sum_{i=1}^M\sum_{n=1}^{N}\sum_{\ell=0}^{n-1}\sum_{m=0}^\ell\kappa_n p_{i,j}  \frac{\mu_j^\ell }{\ell!}\binom{\ell}{m} (-1)^\ell\bigl(\lta_i(\mu_j)\bigr)^{(\ell-m)} \bigl(\ltw_i(\mu_j)\bigr)^{(m)}
\end{equation}
and the waiting-time density is given by
$$
f_{W,j}(x)=\sum_{i=1}^M  \sum_{n=1}^{N}\sum_{\ell=0}^{n-1}\sum_{m=0}^\ell\kappa_n  p_{i,j} \frac{(-1)^\ell}{\ell!}\binom{\ell}{m}\bigl(\lta_i(\mu_j)\bigr)^{(\ell-m)}  \bigl(\ltw_i(\mu_j)\bigr)^{(m)}\mu_j^n \mathrm{e}^{-\mu_j x} \frac{x^{n-\ell-1}}{(n-\ell-1)!}.
$$
\end{theorem}

\begin{proof}
For the proof, we shall refrain from presenting detailed computations, as the analysis is straightforward and similar to the one for the exponential case. We give, however, a few intermediate formulas. From \eqref{7eq:ltw e3is} and for the present preparation time distributions we have that
\begin{multline*}
\ltw[n+1]_{j}(s)=\sum_{i=1}^M \p[Z_n=i] p_{i,j}\Biggl(1 -\sum_{n=1}^{N}\kappa_n\e[\mathrm{e}^{-\mu_j(A_n+W_n)}\sum_{\ell=0}^{n-1}\frac{\mu_j^\ell (A_n+W_n)^\ell}{\ell !}\mid Z_n=i]+\\
+\sum_{n=1}^{N} \kappa_n \Bigl(\frac{\mu_j}{\mu_j+s}\Bigr)^{n}\e[\mathrm{e}^{-\mu_j(A_n+W_n)} \sum_{\ell=0}^{n-1}\frac{(\mu_j+s)^\ell (A_n+W_n)^\ell}{\ell !}\mid Z_n=i]\Biggr).
\end{multline*}
So for $n\to\infty$ we have that $\ltw_j(s)$ is given by (cf.\ \eqref{7eq:ltw B exp})
\begin{equation}\label{7eq:ltw B PH}
\ltw_j(s)=u_j+\sum_{i=1}^M  \sum_{n=1}^{N}\sum_{\ell=0}^{n-1}\sum_{m=0}^\ell\kappa_n  p_{i,j} \frac{(-\mu_j)^\ell}{\ell!}\binom{\ell}{m}\bigl(\lta_i(\mu_j)\bigr)^{(\ell-m)}  \bigl(\ltw_i(\mu_j)\bigr)^{(m)}\Bigl(\frac{\mu_j}{\mu_j+s}\Bigr)^{n-\ell},
\end{equation}
where $u_j$ is the mass of the waiting-time distribution at the origin and it is given by \eqref{mass}.
Inverting the Laplace transform $\ltw_j$ yields that the density of the waiting time is given by the expression presented in the theorem and the corresponding distribution has mass $u_j$ at the origin. For $i,j=1,\ldots,M$, the $M^2$ unknown constants $\ltw_{i}(\mu_j)$ are the unique solution to the system of linear equations resulting from substituting $\mu_k$, $k=1,\ldots,M$, for $s$  in \eqref{7eq:ltw B PH}.
\end{proof}

Naturally, Theorem~\ref{7th:MC-PH B} reduces to Theorem~\ref{7th:MC-exp B} for $N=1$, and to the i.i.d.\ case for $M=1$; cf.\ \cite[Theorem~2]{vlasiou05}.

Dependence structures of the form of Equation~\eqref{7eq:2}, and several generalisations, have been studied extensively for Lindley's recursion. Asmussen and Kella~\cite{asmussen00} considered multidimensional martingales for a class of processes, and sketched the results for the Markov-modulated M/G/1 queue. The basic model was then described in detail in Adan and Kulkarni~\cite{adan03}, where all results are given explicitly and the analysis extends to the study of the queue length distribution. Specifically, the authors study a single-server queue where the in\-ter\-ar\-ri\-val times and the service times depend on a common discrete-time Markov chain in a way similar to the one described by Equation~\eqref{7eq:2}. This model generalises the well-known MAP/G/1 queue by allowing dependencies between interarrival and service times. The MAP/G/1 queue provides a powerful framework to model dependencies between successive interarrival times \cite{combe98}, but typically the service times are assumed to be i.i.d.\ and independent of the arrival process.

\begin{remark}[Periodic Markov chain]
The results in Sections \ref{7ss:Exponential preparation times} and \ref{7ss:Phase-type preparation times} are also valid in case the transition probability matrix $\mathbf{P}$ is periodic. Clearly, the transforms $\ltw_{j}$ should then be defined as Cesaro limits, i.e., for $\mathrm{Re}(s)\geqslant 0$, $n\geqslant 1$, and $j=1,2,\ldots,M$,
$$
\ltw_{j}(s)=\lim_{n\to\infty} \frac{1}{n} \sum_{m = 1}^n \ltw[m]_{j}(s).
$$
But contrary to the aperiodic case, the autocorrelation between service times and preparation times does not converge to zero as the lag tends to infinity; see (\ref{eq:autoA}).
\end{remark}

\subsection{Numerical examples}\label{ss:egs1}
In this section we present some examples to demonstrate the effects of autocorrelation and cross\-cor\-relation of the preparation and service times. In all of the examples, we have assumed that the underlying Markov chain has four states, i.e., $M=4$. The transition matrices are chosen so that the limiting distributions for all cases are equal. Naturally, different transition matrices produce different dependence structures between and among the sequences $\{A_n\}$ and $\{B_n\}$. We also assume that the service times and the preparation times are exponentially distributed and that the service rates in each state of the Markov chain are given by $(\lambda_1, \lambda_2, \lambda_3, \lambda_4)=(1, 100, 1, 100)$. We keep the mix of small and large preparation and service rates the same. Thus the examples differ only in the dependence structure among and between the sequences of preparation and service times.\\

\begin{figure}
\begin{center}
\includegraphics[width=0.8\textwidth]{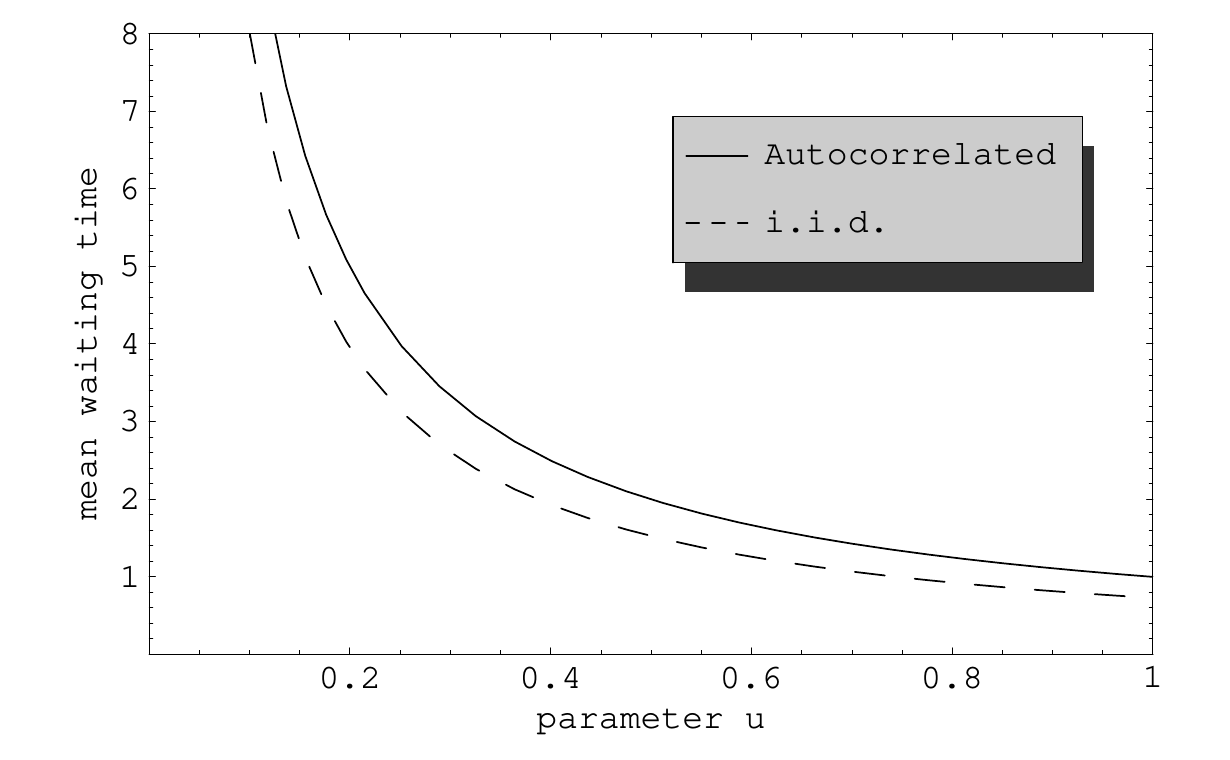}
\caption{Positive crosscorrelation: the mean waiting time against the parameter controlling the preparation rate for autocorrelated and i.i.d.\ sequences.}
\label{fig:1}
\end{center}
\end{figure}

\noindent
{\textbf{First example.} In this example we study the effect of autocorrelation in case $\{A_n\}$ and $\{B_n\}$ are positively crosscorrelated. In particular, we assume that the preparation rates are given by $(\mu_1, \mu_2, \mu_3, \mu_4)=(u/2, 10 u, u/2, 10 u)$, where we employ the parameter $u>0$ in order to explore the effect of the mean preparation time on the expected waiting time. For this setting, we have that the crosscorrelation between $\{A_n\}$ and $\{B_n\}$ is approximately equal to 0.3195.

We will compare the case where both the preparation and the service times are autocorrelated to the case that successive preparation and service times are i.i.d. In the first case, we take the transition probability matrix to be given by
\begin{equation}\label{mat:1}
\mathbf{P}=\begin{pmatrix}
  0 &1 &0 &0\\
  0 &0 &1 &0\\
  0 &0 &0 &1\\
  1 &0 &0 &0\\
\end{pmatrix}
\end{equation}
Thus, since we assumed that the service times are exponentially distributed (with the given rates), we have that their autocorrelation is given by
\begin{equation}\label{eq:auto}
\rho[A_1,A_{n+1}]=(-1)^n \frac{9801}{29803},\qquad n\geqslant 1,
\end{equation}
and similarly, we can derive that the autocorrelation for the preparation times is equal to
$$
\rho[B_1,B_{n+1}]=(-1)^n \frac{361}{1163},\qquad n\geqslant 1.
$$
\begin{figure}
\begin{center}
\includegraphics[width=0.8\textwidth]{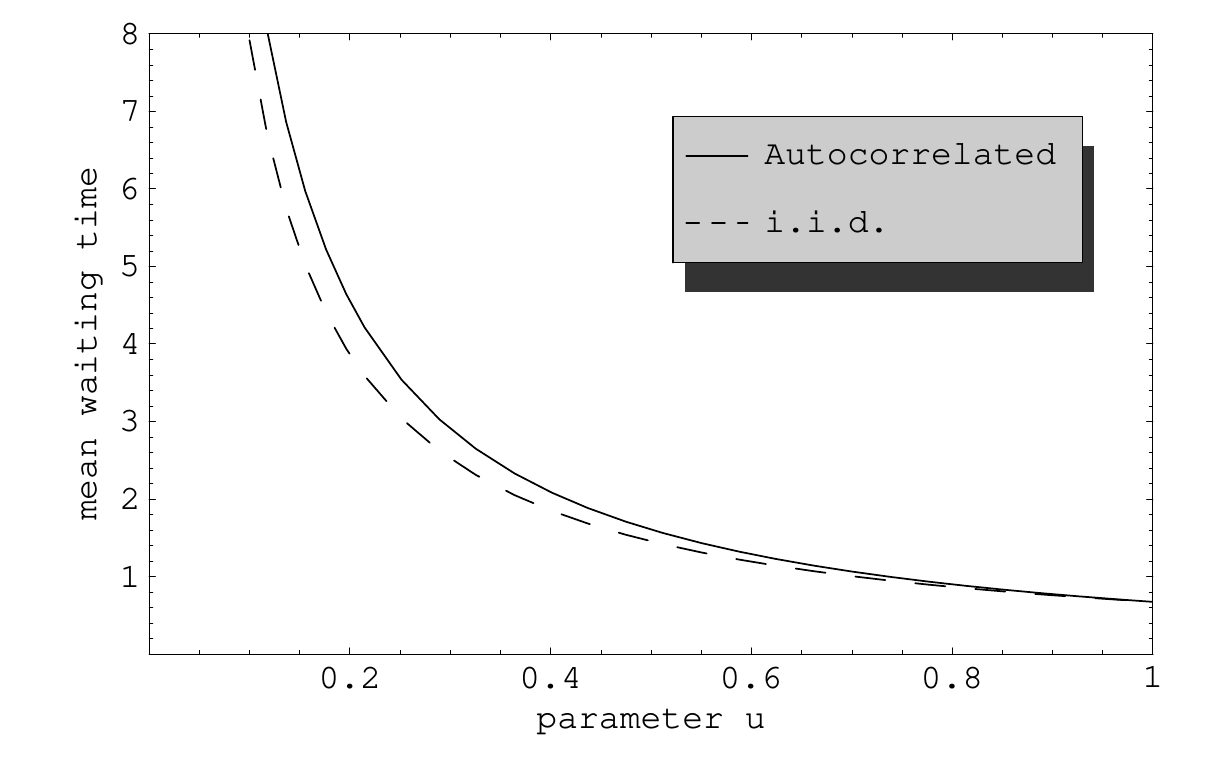}
\caption{Negative crosscorrelation: the mean waiting time against the parameter controlling the preparation rate for autocorrelated and i.i.d.\ sequences.}
\label{fig:2}
\end{center}
\end{figure}
In the second case, we take $\mathbf{P}=(0.25)$, i.e.\ we take all elements of the matrix $\mathbf{P}$ equal to 0.25, which implies that successive preparation and successive service times have zero autocorrelation. In Figure~\ref{fig:1} we plot the mean waiting time of the server against the parameter $u$. Keep in mind that all means of the preparation times decrease as $u$ increases. As is evident from Figure~\ref{fig:1}, autocorrelation leads to higher mean waiting times.\\

\noindent
{\textbf{Second example.} In Figure~\ref{fig:2} we study the effect of autocorrelation in case $\{A_n\}$ and $\{B_n\}$ are \emph{negatively} crosscorrelated. For this case, we take the preparation rates equal to $(\mu_1, \mu_2, \mu_3, \mu_4)=(10 u, u/2, 10 u, u/2)$, $u>0$, which implies that the crosscorrelation in this example is approximately equal to $-0.3195$. We compare the previous two cases for the new preparation rates. Namely, we compare the autocorrelated case, where the transition matrix $\mathbf{P}$ is given by \eqref{mat:1}, to the i.i.d.\ case, where the transition matrix has all its entries equal to 0.25. The autocorrelation functions, both for the preparation times and for the service times remain the same. As before, autocorrelation leads to higher mean waiting times. However, in this case, the effect of the autocorrelation is not as big as in the first example.\\

\begin{figure}[bt]
\begin{center}
\includegraphics[width=0.8\textwidth]{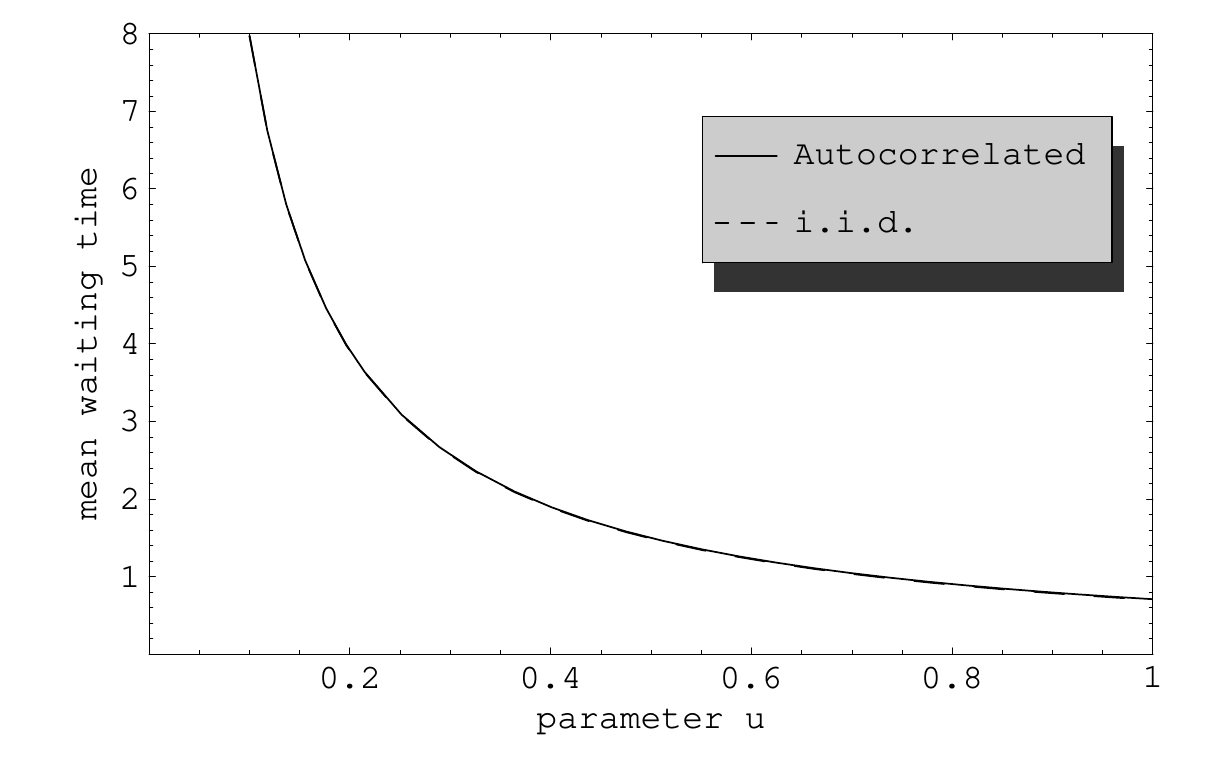}
\caption{Zero crosscorrelation: the mean waiting time against the parameter controlling the preparation rate for autocorrelated and i.i.d.\ sequences.}
\label{fig:3}
\end{center}
\end{figure}

\noindent
{\textbf{Third example.} Now the preparation rates are given by $(\mu_1, \mu_2, \mu_3, \mu_4)=(u/2, u/2, 10 u, 10 u)$, $u>0$, which implies that the crosscorrelation in this example is equal to zero. As before, we compare two cases. In the first case, we take
\begin{equation*}
\mathbf{P}=\begin{pmatrix}
  0 &0.5 &0 &0.5\\
  0.5 &0 &0.5 &0\\
  0 &0.5 &0 &0.5\\
  0.5 &0 &0.5 &0\\
\end{pmatrix}
\end{equation*}
which leads to autocorrelated service times with the autocorrelation function given by Equation~\eqref{eq:auto}, and i.i.d.\ preparation times. We compare this case to the case where all entries of the transition matrix are equal to 0.25, where successive preparation and successive service times have zero autocorrelation. In Figure~\ref{fig:3}, the mean waiting times are shown for both cases. The difference in the mean waiting times is negligible in this case, although the two functions are not equal.

Comparing these numerical results to the analogous cases for Lindley's recursion, which are described in \cite{adan03}, illustrates the different effect of dependencies on these two recursions. In Lindley's recursion, autocorrelation leads to slightly lower mean waiting times. Moreover, negative crosscorrelation leads to a bigger difference between the autocorrelated and the i.i.d.\ case than the difference observed for the positively crosscorrelated case. The results for the Lindley-type recursion~\eqref{recursion} are reversed. Autocorrelation leads to higher mean waiting times and the biggest differences are observed in case $\{A_n\}$ and $\{B_n\}$ are positively crosscorrelated (cf.\ Figures~\ref{fig:1} and \ref{fig:2}). Also, in the third example, where the crosscorrelation is eliminated, the outcome is again very different. For Lindley's recursion, even in the case of zero crosscorrelation, autocorrelation leads to significant differences for the mean waiting times. In our case, there is hardly any difference between the autocorrelated and the i.i.d.\ case (cf.\ Figure~\ref{fig:3}). Keep in mind that contrary to the single server queue, the model studied here assumes an infinite number of customers. Thus, for example, $A_n$ represents the $n$-th service time rather than the $n$-th interarrival time.

\section{Services depending on the previous preparation time}\label{7s:B_n+1 depends on A_n}
In this section, we study the second dependence structure mentioned in the introduction. In particular, we assume that for all $n$, the service times $A_n$ are distributed as $A$, which in turn is exponentially distributed with rate $\lambda$. Moreover, the Laplace-Stieltjes transform of the preparation time $B_{n+1}$, given that the previous service time $A_n$ equals $t$, is of the form
\begin{equation}\label{7eq:dependence structure}
\e[\mathrm{e}^{-s B_{n+1}}\mid A_n=t]=\e[\mathrm{e}^{-s B}\mid A=t]=\chi(s)\mathrm{e}^{-\psi(s) t}.
\end{equation}
Observe that now the preparation time $B_{n+1}$ depends only on the previous service time, while in the Markov-modulated case we have examined previously all preparation and service times are correlated between and among one another, since their distributions depend on a common Markov chain.

This dependence structure also occurs in simple queuing models. Consider the following situation. Work arrives at a single-server queue according to a process with stationary, non-negative independent increments. This work, however, does not immediately enter the queue of the server facility; instead it is accumulated behind a gate. At exponential interarrivals the gate is opened and -- after the addition of an independent component -- the work is collected and delivered as a single customer at the queue of the service facility. The additional component may be viewed as a set-up time.

Due to the exponentially distributed interarrival times of customers, we can view the service facility as an M/G/1 queue in which the interarrival and service time for each customer are positively correlated. Indeed, if the interval between two consecutive openings of the gate is relatively long (short), it is likely that a relatively large (small) amount of work has accumulated during that interval. This model is a unification and generalisation of the M/G/1 queue with a positive correlation between interarrival and service times~\cite{borst92,borst93,cidon96,conolly79} and has been analysed in Boxma and Comb\'{e}~\cite{boxma93}. In Comb\'{e} and Boxma~\cite{combe98} it is shown that the collect system can also be modelled by using the BMAP framework.

\subsection{The waiting-time distribution}
For this model Vlasiou~\cite{vlasiou07a} established (under the most general conditions) that there exists a unique steady-state waiting-time distribution. In order to derive the steady-state waiting-time distribution we shall further assume that the functions $\chi$ and $\psi$ that appear in Equation~\eqref{7eq:dependence structure} are rational functions; i.e.,
\begin{equation}\label{7eq:chi psi}
\chi(s)=\frac{P_1(s)}{Q_1(s)}\qquad\mbox{and}\qquad\psi(s)=\frac{P_2(s)}{Q_2(s)},
\end{equation}
where $P_2$, $Q_1$, and $Q_2$ are polynomials of degrees $L$, $M$, and $N$ respectively. From the form of Equation~\eqref{7eq:dependence structure} we see that a number of other assumptions have been implicitly made. For example, since for $s=0$ the expectation $\e[\mathrm{e}^{-s B}\mid A=t]$ should be equal to one, we have implicitly assumed that $\psi(0)=0$ and $\chi(0)=1$. We shall mention other implications of such type only when necessary.

The preparation time $B_{n+1}$ consists of two parts: a component which depends on the previous service time, represented by $\mathrm{e}^{-\psi(s) t}$, and an `ordinary' preparation time with Laplace-Stieltjes transform $\chi(s)$, which does not depend on the service time. From \eqref{7eq:dependence structure} we have that the bivariate Laplace-Stieltjes transform of the generic preparation and service time is given by
\begin{equation}\label{7eq:LTS B minus A}
\e[\mathrm{e}^{-s B-z A}]=\int_0^\infty \lambda {\rm e}^{-\lambda t} \mathrm{e}^{-z t} \chi(s) \mathrm{e}^{-\psi(s) t}\,{\rm d}t=\frac{\lambda  \chi(s)}{\lambda +\psi(s)+z},
\end{equation}
for $\mathrm{Re}(\lambda +\psi(s)+z)>0$. This expression leads to
$$
\e[B]=\frac{\psi^\prime(0)}{\lambda }-\chi^\prime(0),\quad\mbox{and}\quad \e[AB]=\frac{2\psi^\prime(0)-\lambda \chi^\prime(0)}{\lambda^2},
$$
from which we have that the covariance between a preparation time and a service time is given by
$$
\cov[A,B]=\frac{\psi^\prime(0)}{\lambda^2}.
$$
The correlation between these two variables can be also computed, see Boxma and Comb\'{e}~\cite{boxma93}. Thus, one can construct a distribution function $\Dfb$ that has any covariance (or correlation) between $A$ and $B$ that is desired. The waiting-time density for this case is given by the following theorem.

\begin{theorem}
Let $F_A (x) = 1 - e^{-\lambda x}$. Under the assumptions \eqref{7eq:dependence structure} and \eqref{7eq:chi psi}, we have that
\[
F_W (0) = c_0
\]
and the waiting-time density is given by
\begin{equation*}
\dfw(x)=\sum_{i=1}^{K} c_i \mathrm{e}^{r_i x}.
\end{equation*}
In the expression above, the constants $r_i$ are the $K$ zeros of the equation
$$
Q_1(s)\bigl((\lambda -s)Q_2(s)+P_2(s)\bigr)=0
$$
and the coefficients $c_i$ are derived explicitly as the solution to a linear system of equations.
\end{theorem}

\begin{proof}
In order to derive the steady-state waiting-time distribution, we shall first derive the Laplace-Stieltjes transform of $\Dfw$. We follow a method based on Wiener-Hopf decomposition. A straightforward calculation, starting from the steady-state version of (\ref{recursion}),
\[
 W\stackrel{\m{D}}{=}\max\{0, B-A-W\},
\]
yields that
\begin{align*}
\ltw(s)&=\e[{\rm e}^{-s W}]=\p[B\leqslant W+A]+\e[\mathrm{e}^{-s(B-W-A)}]-\e[\mathrm{e}^{-s(B-W-A)} \,; B\leqslant W+A]\\
        &=\p[B\leqslant W+A]+\e[\mathrm{e}^{s W}]\,\e[\mathrm{e}^{-s(B-A)}]-\e[\mathrm{e}^{-s(B-W-A)} \,; B\leqslant W+A],
\end{align*}
since $B-A$ and $W$ are independent. Therefore, from \eqref{7eq:LTS B minus A} we have for $\mathrm{Re}(s)=0$ that
\begin{align}
\label{7eq:1}  \ltw(s)&=\p[B\leqslant W+A]+\ltw(-s)\frac{\lambda  \chi(s)}{\lambda -s+\psi(s)}-\e[\mathrm{e}^{-s(B-W-A)} \mid B\leqslant W+A]\p[B\leqslant W+A]\\
\nonumber        &=\ltw(-s)\frac{P_1(s)}{Q_1(s)}\frac{\lambda Q_2(s)}{(\lambda -s)Q_2(s)+P_2(s)}+\p[B\leqslant W+A]\bigl(1-\e[\mathrm{e}^{-s(B-W-A)} \mid B\leqslant W+A]\bigr),
\end{align}
which can be rewritten as
\begin{multline}\label{7eq:ltw etoimo gia WH}
\ltw(s)Q_1(s)\bigl((\lambda -s)Q_2(s)+P_2(s)\bigr)=
\lambda \ltw(-s)P_1(s) Q_2(s)+Q_1(s)\bigl((\lambda -s)Q_2(s)+P_2(s)\bigr)\times\\
\times\p[B\leqslant W+A]\bigl(1-\e[\mathrm{e}^{-s(B-W-A)} \mid B\leqslant W+A]\bigr).
\end{multline}
We can observe that $Q_1(s)\bigl((\lambda -s)Q_2(s)+P_2(s)\bigr)$ is a polynomial of degree $K=\max\{M+N+1, M+L\}$ and also that the left-hand side of \eqref{7eq:ltw etoimo gia WH} is analytic for $\mathrm{Re}(s)>0$ and continuous for $\mathrm{Re}(s)\geqslant0$, and the right-hand side of \eqref{7eq:ltw etoimo gia WH} is analytic for $\mathrm{Re}(s)<0$ and continuous for $\mathrm{Re}(s)\leqslant0$. Therefore, from Liouville's theorem \cite{titchmarsh-tf}, we conclude that both sides of \eqref{7eq:ltw etoimo gia WH} are the same $K$-th degree polynomial, say, $\sum_{i=0}^{K} q_i s^i$. Hence,
\begin{equation}\label{7eq:transform in poly}
  \ltw(s)=\frac{\sum_{i=0}^{K} q_i s^i}{Q_1(s)\bigl((\lambda -s)Q_2(s)+P_2(s)\bigr)}.
\end{equation}
In the expression above, the constants $q_i$ are not determined so far. In order to obtain the transform, observe that $\ltw$ is a fraction of two polynomials both of degree $K$. Let $r_i$, $i=1,\ldots,K$, be the zeros of the denominator. Ignoring the special case of zeros with multiplicity greater than one, partial fraction decomposition yields that \eqref{7eq:transform in poly} can be rewritten as
\begin{equation}\label{7eq:transform in fractions}
  \ltw(s)=c_0+ \sum_{i=1}^{K} \frac{c_i}{s-r_i},
\end{equation}
which implies that the waiting-time distribution has a mass at the origin that is given by
$$
\p[W=0]=\lim_{s\to\infty}\e[\mathrm{e}^{-sW}]=c_0
$$
and has a density that is given by
\begin{equation}\label{7eq:density}
\dfw(x)=\sum_{i=1}^{K} c_i \mathrm{e}^{r_i x}.
\end{equation}

All that remains is to determine the $K+1$ constants $c_i$. To do so, we work as follows. We express the terms $\p[B\leqslant W+A]$ and $\e[\mathrm{e}^{-s(B-W-A)} \mid B\leqslant W+A]$ that appear at the right-hand side of \eqref{7eq:ltw etoimo gia WH} in terms of the constants $c_i$ by using \eqref{7eq:transform in fractions} and \eqref{7eq:density}. Then, we substitute these expressions and \eqref{7eq:transform in fractions} in the left-hand side of \eqref{7eq:ltw etoimo gia WH}. Thus we obtain a new equation in terms of the constants $c_i$ that we differentiate a total of $K$ times. We evaluate each of these derivatives for $s=0$ and thus we obtain a linear system of $K+1$ equations, $i=0,\ldots,K$, for the constants $c_i$ (the original equation and the $K$ derivatives).
\end{proof}

\begin{remark}
In case that for some $i$ we have that the real part of $r_i$ is greater than or equal to zero, then it follows that the corresponding coefficient $c_i$ is equal to zero. This is a consequence of the fact that the limiting distribution exists and therefore the density $\dfw$ should be integrable. Note that the meaning of $\chi$ as the Laplace-Stieltjes transform of an ordinary preparation time already implies that the $L$ zeroes of $Q_1(s)$ have a negative real part.
\end{remark}
\begin{remark}\label{7rem:complex roots}
Although the roots $r_i$ and coefficients $c_i$ may be complex, the density and the mass $c_0$ at zero will be positive. This follows from the fact that there is a unique equilibrium distribution and thus a unique solution to the linear system for the coefficients $c_i$. Of course, it is also clear that each root $r_i$ and coefficient $c_i$ have a companion conjugate root and conjugate coefficient, which implies that the imaginary parts appearing in the density cancel.
\end{remark}
\begin{remark}
When $Q_1(s)\bigl((\lambda -s)Q_2(s)+P_2(s)\bigr)$ has multiple zeros, the analysis proceeds in essentially the same way. For example, if  ${r}_1={r}_2$, then the partial fraction decomposition of $\ltw$ becomes
$$
\ltw(s)=c_0+ \frac{c_1}{(s-r_1)^2}+\sum_{i=2}^{K} \frac{c_i}{s-r_i},
$$
yielding
$$
\dfw(x)={c}_1 x \mathrm{e}^{r_1 x}+\sum_{i=2}^{K} c_i \mathrm{e}^{r_i x}.
$$
\end{remark}
\begin{remark}
For the service time $A$ we have considered only the exponential distribution, mainly because we can illustrate the technique we use without complicating the analysis. However, we can extend this class by considering distributions with a mixed-Erlang distribution of the form of Equation~\eqref{7eq:distr} and the proof remains essentially the same. That is, let $\Dfa(x)=\sum_{i=1}^n \kappa_i E_i(x)$, where $E_i$ is the Erlang distribution with $i$ phases. Then, the resulting density of the waiting time is again of the form
$$
\dfw(x)=\sum_{i=1}^{K^\prime} c_i \mathrm{e}^{r_i x},
$$
where $K^\prime=M+n\max\{N+1, L\}$. The constants $r_i$ are the zeros to the equation
$$
Q_1(s)\bigl((\lambda -s)Q_2(s)+P_2(s)\bigr)^n=0,
$$
and the coefficients $c_i$ are determined in a similar fashion as before. Naturally, if any of the zeros $r_i$ has multiplicity greater than one, the form of the density changes analogously; see also the previous remark.
\end{remark}
We now present a few examples where we show how some classic dependence structures fit into this class. In the examples below, we only discuss how to derive the function $\psi$. All of the examples presented here are examples of L\'{e}vy processes.

\pdfbookmark[2]{Independent case}{anexarthto}
\subsection*{Independent case}
The dependence structure described by Equation~\eqref{7eq:dependence structure} contains a great variety of dependence structures, including the independent case. If for all $s$ we have that $\psi(s)=0$, then the Laplace-Stieltjes transform of $B$ is independent of the length of the service time, and thus the function $\chi$ appearing in \eqref{7eq:dependence structure} is in fact the Laplace-Stieltjes transform of $B$. Observe that we have assumed that $B$ has a rational Laplace-Stieltjes transform, which is necessary in order to decompose Equation~\eqref{7eq:1} into functions that are analytic either in the left-half plane or in the right-half plane.

\pdfbookmark[2]{Linear Dependence}{grammiko}
\subsection*{Linear Dependence}
Assume that the service time $A$ and the preparation time $B$ are linearly dependent; that is, $B=c A$. Then,
$$
\e[\mathrm{e}^{-s B}\mid A=t]=\e[\mathrm{e}^{-s cA}\mid A=t]=\mathrm{e}^{-s ct}.
$$
Thus we have that $\chi(s)=1$, and $\psi(s)=cs$, and both functions satisfy our assumptions. Queuing models with a linear dependence between the service time and the preceding interarrival time have been studied in \cite{cidon96,conolly68,conolly69}; see also \cite{boxma01}.
\pdfbookmark[2]{Compound Poisson Process}{mikthpsari}
\subsection*{Compound Poisson Process}
In this case we assume that given that $A=t$, the preparation time $B$ is equal to $\sum_{i=1}^{N(t)} C_i$, where $N(t)$ is a Poisson process with rate $\gamma$, and $\{C_i\}$ is a sequence of i.i.d.\ random variables, where each of them is distributed like $C$, and where $C$ has a rational Laplace-Stieltjes transform. By convention, if $N(t)=0$, we have that the preparation time is also zero. Under this assumption, we have
\begin{align*}
\e[\mathrm{e}^{-s B}\mid A=t]&=\e[\mathrm{e}^{-s \sum_{i=1}^{N(t)} C_i}\mid A=t]
                            =\sum_{k=0}^\infty \e[\mathrm{e}^{-s \sum_{i=1}^{k} C_i}\mid A=t]\,\mathrm{e}^{-\gamma t}\frac{(\gamma t)^k}{k!}\\
                            &=\sum_{k=0}^\infty \bigl(\e[\mathrm{e}^{-s C}]\bigr)^k\,\mathrm{e}^{-\gamma t}\frac{(\gamma t)^k}{k!}=\mathrm{e}^{-\psi(s) t},
\end{align*}
where $\psi(s)=\gamma \bigl(1-\e[\mathrm{e}^{-s C}]\bigr)$. As before, in this case we have that $\chi(s)=1$.

\pdfbookmark[2]{Brownian Motion}{kafes}
\subsection*{Brownian Motion}
In this case we assume that given that $A=t$, the random variable $B$ is normally distributed with mean $\mu t$ and variance $\sigma^2 t$. Then we have that
$$
\e[\mathrm{e}^{-s B}\mid A=t]=\int_{-\infty}^\infty \mathrm{e}^{-s x}\frac{\mathrm{e}^{-{(x-\mu t)^2}/{(2\sigma^2 t)}}}{\sigma \sqrt{2 \pi t}}\,{\rm d}x=\mathrm{e}^{-\psi(s) t},
$$
where $\psi(s)=\mu s-{s^2\sigma^2}/{2}$, and $\chi(s)=1$. Of course, if $B$ is interpreted as the preparation time of a customer, then assuming that $B$ is normally distributed is not a natural assumption, since the preparation time of a customer is non-negative. However, in the analysis we do not need the condition of $B$ being non-negative; therefore, it is mathematically possible to consider this case. For further examples of distributions satisfying the condition described by \eqref{7eq:dependence structure}, see Boxma and Comb\'{e}~\cite{boxma93}.

One additional assumption made in \cite{boxma93} is that the function $\psi$ appearing in \eqref{7eq:dependence structure} has a completely monotone increasing derivative. This implies that the function $\mathrm{e}^{-\psi(s)}$ is the Laplace-Stieltjes transform of an infinitely divisible probability distribution, which in its turn is a distribution of increments in processes with stationary independent increments, i.e.\ L\'{e}vy processes. However, the monotonicity of the derivative of $\psi$ was not required for the analysis of our model.

\subsection{Numerical examples}\label{ss:egs2}
We shall demonstrate the effect of the dependence structure given by \eqref{7eq:dependence structure} to the mean waiting time of the server. In Figure~\ref{fig2} we plot the mean waiting time against the mean preparation time of a customer for the three first cases mentioned in the previous section. In all cases, we assume that the service times are exponentially distributed with mean 1.

As an example, we present the derivations for the linear dependence. In this case, we have that $K=1$ and thus we have one root, given by $r_1=-1/(c-1)$. Moreover, one can compute that
$$
\p[B\leqslant W+A]=\begin{cases}
  1, &c\leqslant 1\\
  1-c_0+\frac{c_1}{2 r_1}, &c> 1
\end{cases}
$$
and that
$$
\e[\mathrm{e}^{-s(B-W-A)} \mid B\leqslant W+A]=\begin{cases}
  \frac{\ltw(-s)}{1-(1-c)s}, &c\leqslant 1\\
  \frac{1}{1+(c-1)s} \frac{\ltw(-s)-\ltw(1/(c-1))}{1-c_0+\frac{c_1}{2 r_1}},&c> 1.
\end{cases}
$$

\begin{figure}[ht]
\begin{center}
\includegraphics[width=0.8\textwidth]{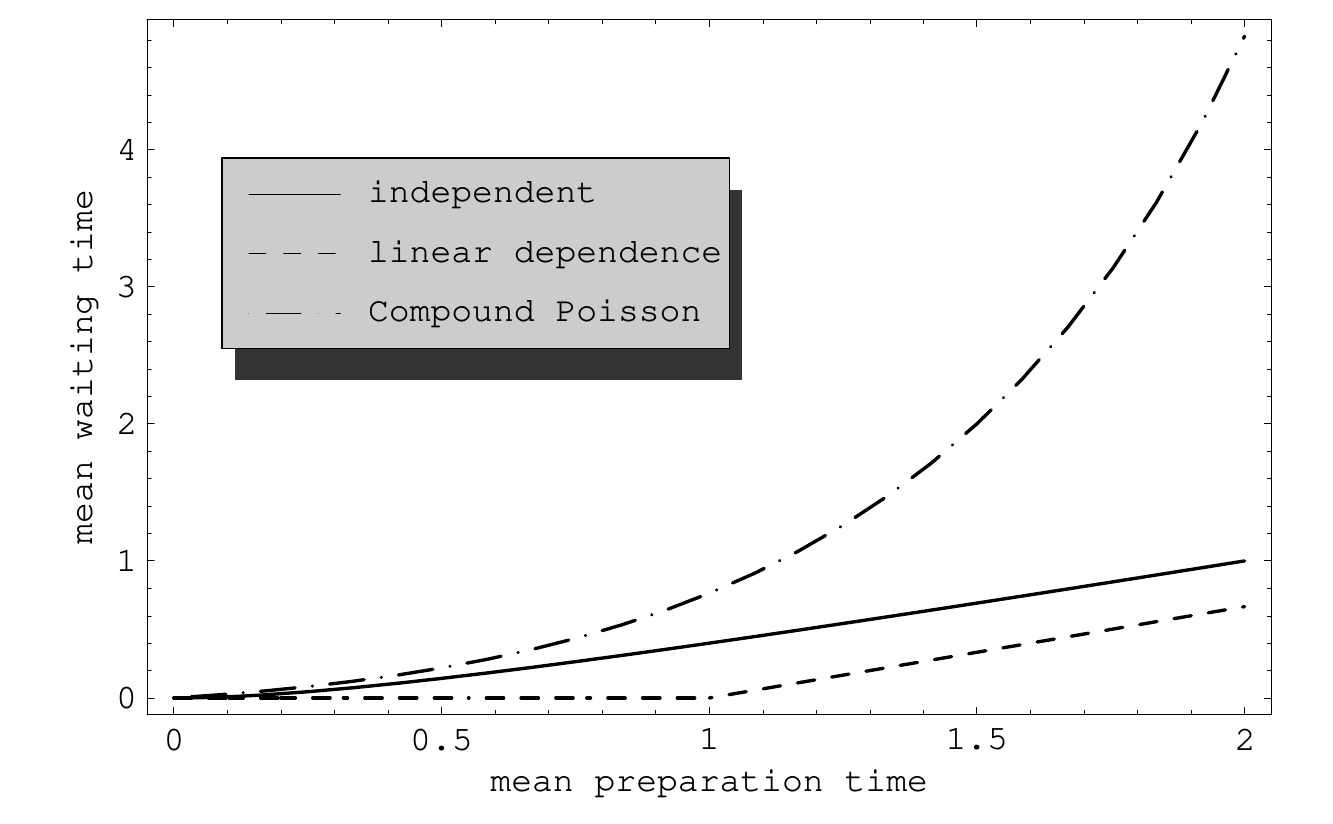}
\caption{The mean waiting time against the mean preparation time}
\label{fig2}
\end{center}
\end{figure}

Thus, for $c\leqslant 1$, Equation~\eqref{7eq:ltw etoimo gia WH} becomes
$$
(c_0+\frac{c_1}{s-r_1})(1 -s+c s)=
c_0-\frac{c_1}{s+r_1}+(1 -s+c s)\bigl(1-\frac{c_0-\frac{c_1}{s+r_1}}{1-(1-c)s}\bigr).
$$
Evaluating this equation and its derivative at $s=0$ we form a linear system of two equations, from which we conclude that for $c\leqslant 1$, $c_0=1$ and  $c_1=0$. This is what was to be expected, as in this case there is no waiting time.

For $c>1$, Equation~\eqref{7eq:ltw etoimo gia WH} becomes
$$
(c_0+\frac{c_1}{s-r_1})(1-s+cs)=c_0-\frac{c_1}{s+r_1}+(1 -s+c s)(1-c_0+\frac{c_1}{2 r_1})\bigl(1- \frac{\frac{c_1}{2 r_1}-\frac{c_1}{s+r_1}}{(1+(c-1)s)(1-c_0+\frac{c_1}{2 r_1})}\bigr).
$$
Evaluating this equation and its derivative at $s=0$ we have that for $c>1$, $c_0=1/3$ and  $c_1=2/(3(c-1))$. Since the mean waiting time is equal to $c_1/r_1^2$, we have that $\e[W]$ is zero for $c\leqslant 1$ and is equal to $2(c-1)/3$ for $c>1$. This case is depicted in Figure~\ref{fig2}.

Furthermore, for the independent case we have assumed that the preparation times are also exponentially distributed, and for the compound Poisson process we have assumed that the jumps $C_i$ are exponentially distributed with mean 1 (which implies that the mean preparation times are equal to the parameter $\gamma$ of the underlying Poisson process).

As is evident from Figure~\ref{fig2}, the precise nature of the dependence structure between the preparation and the service times has a significant impact on the mean waiting time (and on the waiting-time distribution). In Section~\ref{ss:egs1}, where we have examined the effect of Markov modulated dependencies, correlation has led to longer waiting times. However, in this case, we observe that dependent preparation and service times do not necessarily lead to longer waiting times. In case $B$ is linearly dependent on $A$, the resulting mean waiting time is always smaller than the mean waiting time for $A$ and $B$ being both exponentially distributed and independent of one another.

We further observe that the mean waiting time can vary significantly among various dependence structures. For example, we see that for the compound Poisson case, the mean waiting time grows faster than an exponential (which is the growth rate for the independent case) and diverges from the independent case faster than the mean waiting time for the linear case. One can infer from this example that the dependence among preparation and service times in Equation~\eqref{recursion} cannot be ignored and should be modelled appropriately.

\phantomsection
\addcontentsline{toc}{section}{References}
\bibliographystyle{abbrv}
\bibliography{maria}

\end{document}